
\documentstyle{amsppt}
\magnification=1200
\input amstex
\NoPageNumbers
\vsize=23.5truecm\hsize=15.5truecm
\hoffset = 0.3truecm

\def\Ext{\operatorname{Ext}} \def\Tor{\operatorname{Tor}} \def\Hom{\operatorname{Hom}} \def\gra{\operatorname{grade}} \def\Coker{\operatorname{Coker}} \def\Ker{\operatorname{Ker}} \def\Im{\operatorname{Im}} \def\depth{\operatorname{depth}} \def\dim{\operatorname{dim}} \def\height{\operatorname{ht}} \def\proj{\operatorname{proj.dim}} \def\Ann{\operatorname{Ann}} \def\rank{\operatorname{rank}}

\topmatter \title \endtitle \thanks{\noindent The authors are partially supported by the National Basic Research Program}\endthanks \endtopmatter

\document \vglue 0.5cm 

\heading \bf SPECIALIZATION OF MODULES \endheading \vglue 0.5cm
\heading Dam Van Nhi \endheading \centerline {Pedagogical College (CDSP), Thai Binh, Vietnam} \smallskip 
\heading Ng\^o Vi\^et Trung \endheading \centerline {Institute of Mathematics, Box 631, B\`o H\^o, Hanoi, Vietnam} \vglue 2.3truecm \baselineskip = 15pt

\heading Introduction \endheading \smallskip

Specialization is a classical method in algebraic geometry. The first ideal-theoretic approach to specialization was introduced by W. Krull in [Kr1], [Kr2]. Given a family $u = (u_1,\ldots,u_m)$ of indeterminates and a family $\alpha = (\alpha_1,\ldots,\alpha_m)$ of elements of an extension of a  base field $k$,  Krull defined the specialization $I_{\alpha}$ of an ideal $I$ of a polynomial ring $R = k(u)[x]$ with respect to the substitution $u \longrightarrow \alpha$ as follows: $$I_{\alpha} =\{f(\alpha,x)|\ f(u,x)\in I\cap k[u,x]\}.$$ This is obviously an ideal of the polynomial ring $R_{\alpha} = k(\alpha)[x]$. For almost all substitutions $u \longrightarrow \alpha$, the ideal $I_{\alpha}$ inherits most of the basic properties of $I$. \smallskip
The notion of specialization of an ideal played an important role in the study of normal varieties by A.Seidenberg [S] who proved that almost all hyperplane sections of a normal variety are normal again under certain conditions, see also W.E. Kuan [Ku1], [Ku2] for local versions of this result.  Following these works, the second author studied the preservation of properties from the quotient ring $R/I$ to $R_{\alpha}/I_{\alpha}$ [T]. It turned out that specializations of ideals behave well in regard to irreducibility, regularity, normality, Cohen-Macaulayness, and Buchsbaumness. However, certain ring-theoretic properties which are characterized by homological means such as Gorensteiness or generalized Cohen-Macaulayness could not be investigated. The reason is that there is no appropriate theory for the specialization of modules which appear in the homological characterizations. \smallskip  
Generalizing the case of ideals, Seidenberg [S] already gave a definition for the specialization of submodules of a free R-module of finite rank. This naturally leads to a definition for the specialization of quotient modules of free $R$-modules of finite ranks and therefore for the specialization of arbitrary finitely generated $R$-modules. But such a definition only makes sense if  it does not depend on the presentation of the given module (up to isomorphisms). For instance, for a submodule of a finite $R$-module of finite rank, one would have two specializations depending on the presentations of the given module as a submodule and as a quotient module of a free module. It is not clear whether these two specializations are isomorphic, and, even when they are isomorphic to each other, the isomorphism should be canonical. \smallskip
The aim of this paper is to give a definition for the specialization of an arbitrary finitely generated $R$-module which, up to canonical isomorphisms, does not depend on the presentation of the given module, and which preserves basic properties and operations on modules.  \smallskip
First, we observe that the specializations of free $R$-modules of finite ranks and of homomorphisms between them can be  defined in a natural way. Then, given an arbitrary finitely generated $R$-module $L$ with a finite  presentation $F_1 \longrightarrow F_0 \longrightarrow L \longrightarrow 0,$ where $F_0, F_1$ are free modules of finite ranks, we define the specialization $L_{\alpha}$ of $L$ to be the cokernel of the map $(F_1)_\alpha \longrightarrow (F_0)_\alpha$. \smallskip
This definition is not as explicit as Seidenberg's definition. However, it has the advantage that all encountered problems concerning the uniqueness of $L_{\alpha}$ (up to isomorphisms) and compatibility with basic operations on modules can be brought back to problems on specializations of free modules and homomorphisms between them which can be studied by the existing theory of specializations of ideals. The uniqueness of $L_{\alpha}$ can be deduced from the preservation of the exactness of  finite complexes of free modules at specializations which, by the criterion of Buchsbaum and Eisenbud [B-E], can be proved by considering the ranks and the depths of determinantal ideals of the given complex. This preservation can be extended to any exact sequence of finitely generated modules. Using this fact we are able to prove that specializations preserve basic operations on modules including the Tor and Ext functors. \smallskip 
For later applications one should also develop a theory for specializations of  finitely generated modules over a local ring $R_P$, where $P$ is an arbitrary prime ideal of $R$. But this is more subtle since the specialization $P_{\alpha}$ of $P$ need not to be a prime ideal in $k(\alpha)[x]$ (see e.g. [T]). We hope that the results of this paper would help us to settle this problem in the near future. \smallskip
This paper is divided into four sections. \smallskip  
In Section 1 we define the specialization of free modules of finite rank and of homomorphisms between them. The main result is the preservation of finite exact complexes of free modules of finite rank (Theorem 1.5). \smallskip  
In Section 2 we introduce the specialization of arbitrary finitely generated $R$-modules and of homomorphisms between them. The specialization of ideals as well as Seidenberg's  specialization of submodules of  free modules of finite ranks can be considered as special cases of our definition. We will prove that specializations of modules preserve the exactness of an exact sequence (Theorem 2.4). \smallskip
In Section 3 we study the preservation of basic operations on modules at specializations. As applications we show that for almost all $\alpha,$ $\dim L_\alpha = \dim L$ (Theorem 3.4) and that primary decomposition specializes if $k$ is an algebraically closed field (Proposition 3.5).\smallskip
In Section 4 we will prove that $$\align \Tor_i^R(L,M)_\alpha &\cong \Tor_i^{R_\alpha}(L_\alpha,M_\alpha),\; i\ge 0,\cr \Ext_R^i(L,M)_\alpha &\cong \Ext_{R_\alpha}^i(L_\alpha,M_\alpha),\; i\ge 0,\endalign$$ \noindent for almost all $\alpha$ (Theorem 4.2 and Theorem 4.3). As a consequence, the grade of an ideal on a module remains unchanged at specializations (Corollary 4.4).\smallskip
Throughout this paper we assume that $k$ is an arbitrary infinite field. We denote by $R$ the polynomial ring $k(u)[x]:= k(u_1,\ldots,u_m)[x_1,\ldots,x_n]$, where $u_1,\ldots,u_m$ are indeterminates. For each $\alpha = (\alpha_1,\ldots,\alpha_m) \in K^m,$ where $K$ is an extension field of $k$, we put $ R_\alpha = k(\alpha)[x]:= k(\alpha_1,\ldots,\alpha_m) [x_1,\ldots,x_n].$  We shall say that a property holds for almost all $\alpha$ if it holds for all $\alpha$ except perhaps those lying on a proper algebraic subvariety of $K^m.$ For other notations we refer the reader to [B-H] or [E].  \medskip

\heading 1. Specialization of free modules \endheading \smallskip

The specialization of a free $R$-module and of a homomorphism between free $R$-modules can be defined in the following natural way.  
\proclaim{Definition} {\rm Let $F$ be a free $R$-module of finite rank. The {\it specialization} $F_\alpha$ of $F$ is a free $R_\alpha$-module of the same rank.}\endproclaim
Let $\phi: F\longrightarrow G$ be a homomorphism of free $R$-modules. Let $\{f_i\}$ and $\{g_j\}$ be  bases of $F$ and $G$, respectively. Write $$\phi(f_i) = \sum_{i=1}^s a_{ij}(u,x)g_j,\; a_{ij}(u,x)\in R.$$ 
\noindent Then we can represent $\phi$ by the matrix $A = (a_{ij}(u,x)).$ \smallskip
Every element $a(u,x) \in R$ can be written in the form $$a(u,x) = {{p(u,x)}\over {q(u)}},\; p(u,x)\in k[u,x],\; q(u)\in k[u]\setminus \{0\}.$$ \noindent For any $\alpha$ such that $q(\alpha)\ne 0$ we define $$a(\alpha,x) = {{p(\alpha,x)}\over {q(\alpha)}}.$$ Set $A_\alpha = (a_{ij}(\alpha,x)).$ Then $A_\alpha$ is well-defined for almost all $\alpha$.  
\proclaim{Definition} \rm{Let $\phi:  F \longrightarrow G$ be a homomorphism of free $R$-modules of finite ranks represented by a matrix $A.$ The {\it specialization} $\phi_\alpha:  F_\alpha \longrightarrow G_\alpha$ of $\phi$ is given by the matrix $A_\alpha$ provided that $A_\alpha$ is well-defined.}\endproclaim
  
\demo{Remark} This definition requires two fixed bases of $F_\alpha$ and $G_\alpha.$ \enddemo \smallskip
If we choose different bases for $F$ and $G$, then $\phi$ is given by a different matrix. However, the above definition of $\phi_\alpha$ does not depend on the choice of the bases of $F$, $G$ in the following sense.
 
\proclaim{ Lemma 1.1} Let  $\phi_\alpha:  F_\alpha \longrightarrow G_\alpha$ be defined as above. Let $B$ be the matrix of $\phi$ with respect to other bases of $F, G$. Then there are bases of $F_\alpha, G_\alpha$ such that $B_\alpha $ is the matrix of $\phi_\alpha$ with respect to these bases for almost all $\alpha$. \endproclaim
\demo{ Proof} There are transformation matrices $C$ and $D$ such that $B = D^{-1}AC$. Since the determinants of $C$ and $D$ are invertible, the determinants of $C_\alpha$ and $D_\alpha$ are invertible for almost all $\alpha.$ Then we may consider $C_\alpha$ and $D_\alpha$ as transformations matrices of $F_\alpha$ and $G_\alpha$, respectively, which map the given bases of $F_\alpha$ and $G_\alpha$ to new bases of $F_\alpha$ and $G_\alpha.$ Clearly, $B_\alpha = D_\alpha^{-1} A_\alpha C_\alpha.$ Hence $B_\alpha$ is the matrix of $\phi_\alpha$ with respect to the new bases of $F_\alpha$ and $G_\alpha.$ \qed \enddemo \smallskip

\proclaim{ Lemma 1.2}  Let $\phi, \psi: F \longrightarrow G$ and $\delta: G \longrightarrow H$ be homomorphisms of free $R$-modules of finite ranks. Then, for almost all $\alpha,$
$$\align (\phi + \psi)_\alpha & = \phi_\alpha + \psi_\alpha, \\ (\delta \phi)_\alpha & = \delta_\alpha \phi_\alpha.\endalign$$ \endproclaim
\demo{Proof} Let  $A,B,C$ be the matrices of $\phi,\psi,\delta$ with respect to fixed bases of $F, G, H$. For almost all $\alpha$ we have $$\align (A+B)_\alpha & = A_\alpha + B_\alpha, \\ (AC)_\alpha & = A_\alpha C_\alpha.\endalign$$ Hence the conclusion is immediate.\qed \enddemo \smallskip

\proclaim{ Lemma 1.3}  Let $\phi, \psi$ be the homomorphisms of free $R$-modules of finite ranks. Then, for almost all $\alpha,$ $$\align (\phi \oplus \psi)_\alpha & = \phi_\alpha \oplus \psi_\alpha,\\ (\phi \otimes \psi)_\alpha & = \phi_\alpha \otimes \psi_\alpha.\endalign $$ \endproclaim  
\demo{Proof} Let $\phi, \psi$ be represented by the matrices $A,B,$ respectively. Then the matrices of $\phi \oplus \psi$ and $\phi \otimes \psi$ are $A \oplus B,$ and $A \otimes B,$ respectively. It is clear that $$\align (A \oplus B)_\alpha & = A_\alpha \oplus B_\alpha,\\  
(A \otimes B)_\alpha & = A_\alpha \otimes B_\alpha \endalign$$ for almost all $\alpha.$\qed \enddemo \smallskip

Let ${\bold F_\bullet}:\; 0 \longrightarrow F_\ell \overset {\phi_\ell} \to {\longrightarrow } F_{\ell-1} \longrightarrow \cdots \longrightarrow F_1 \overset {\phi_1}\to \longrightarrow F_0$ be a finite complex of free $R$-modules finite ranks. By Lemma 1.2, $$({\bold F_\bullet})_\alpha :\; 0 \longrightarrow (F_\ell)_\alpha \overset {(\phi_\ell)_\alpha} \to \longrightarrow  (F_{\ell-1})_\alpha \longrightarrow \cdots \longrightarrow (F_1)_\alpha \overset {(\phi_1)_\alpha}\to \longrightarrow (F_0)_\alpha $$ is a complex of free $R_\alpha$-modules for almost all $\alpha.$ To study the preservation of the exactness of finite complexes of modules at specializations we shall need the following notations. \smallskip
Let $\phi: F \longrightarrow G$ be a homomorphism of free $R$-modules of finite ranks and let $A$ be the matrix of $\phi$ w.r.t. fixed bases of $F,G$. Denote by $I_t (\phi)$ the ideal generated by the $t \times t$ minors of $A$. Put $\rank\phi := \max \{t \vert \; I_t (\phi) \ne 0 \}$ and put $I(\phi):= I_d (\phi)$ if $d = \text{rank}\,\phi.$ Then $ {\bold F_\bullet}$ is exact if and only if $$\align \text{rank}\, F_i &= \rank\phi_i + \rank\phi_{i +1},\\ &\text{\rm depth}\,I(\phi_i)  \ge i \endalign $$ for $i = 1,\ldots, \ell,$ (see [B-E, Corollary 1] or [E, Theorem 20.9]). Therefore,  we need to study $\text{rank}\,\phi_\alpha$ and $\depth I(\phi_\alpha).$ \smallskip

\proclaim{ Lemma 1.4}  Let $\phi:F \longrightarrow G$ be a homomorphism of free $R$-modules of finite ranks. Then, for almost all $\alpha,$ \smallskip
\noindent \rm (i)\quad $\rank\phi_\alpha = \rank\phi,$\smallskip 
\noindent \rm (ii)\quad $\depth I(\phi_\alpha) = \text{\rm depth\,}I(\phi).$ \endproclaim
\demo{ Proof} (i) Put $d$ = rank$\,\phi,$ $d_\alpha$ = rank$\,\phi_\alpha$. Then $d_\alpha \le d$. By the definition of rank$(\phi)$, there is a $d\times d$ minor $D \ne 0$ of $A$. For almost all $\alpha$, $D_\alpha \ne 0.$ Hence $d_\alpha \ge d.$ From this it follows that $d_\alpha = d$.\smallskip
\noindent (ii) Since $R$ and $R_\alpha$ are Cohen-Macaulay rings, $$\align \depth I(\phi)& = \text{ht}\,I(\phi)\\ \depth  I(\phi_\alpha) & = \text{ht}\,I(\phi_\alpha).\endalign$$ By [T, Lemma 1.1], $\text{ht}\,I(\phi) = \text{ht}\,I(\phi)_\alpha$ for almost all $\alpha$. Let $f_1(u,x),\ldots,f_s(u,x)$ be the $d \times d$ minors of $A$. Then $I(\phi) = (f_1(u,x),\ldots,f_s(u,x))$. By [Kr1, \S1, Satz 1], $I(\phi)_\alpha = (f_1 (\alpha,x),\ldots, f_s (\alpha,x))$ for almost all $\alpha$. Since $f_1(\alpha,x),\ldots,f_s(\alpha,x)$ are the $d\times d$ minors of the matrix $A_\alpha$, we also have $I(\phi_\alpha) = (f_1(\alpha,x), \ldots ,f_s(\alpha,x)$.  Thus, $$I(\phi)_\alpha = I(\phi_\alpha)$$ for almost all $\alpha$. Hence $\depth I(\phi_\alpha) = \depth I(\phi)$ for almost all $\alpha$.\qed \enddemo \smallskip

\proclaim {Theorem 1.5} Let $\bold F_\bullet$ be a finite exact complex of free $R$-modules of finite ranks. Then $(\bold F_\bullet)_\alpha$ is an exact complex of free $R_\alpha$-modules of finite ranks for almost all $\alpha.$ \endproclaim 
\demo{ Proof} By Lemma 1.4, we have $$\align \rank&(\text {F}_i)_\alpha  = \rank(\phi_i)_\alpha + \rank(\phi_{i+1})_\alpha,\\ &\depth I((\phi_i)_\alpha)  = \depth I(\phi_i) \ge i,\endalign$$ 
\noindent $i = 1,\ldots,\,\ell,$ for almost all $\alpha.$ Hence $(\bold F_\bullet)_\alpha$ is an exact complex of free $R_\alpha$-modules for almost all $\alpha.$\qed \enddemo \smallskip

\proclaim{Corollary 1.6}  Let $F \overset \phi \to \longrightarrow G \overset \psi \to \longrightarrow H$ be an exact sequence of free $R$-modules of finite ranks. Then $F_\alpha \overset {\phi_\alpha}\to \longrightarrow G_\alpha \overset{\psi_\alpha}\to \longrightarrow H_\alpha$ is an exact sequence for almost all $\alpha$.\endproclaim 
\demo {Proof} Since the base ring $R$ is regular, we may extend the above exact sequence to a finite exact complex of free $R$-modules [E, Corollary 19.8], [B-H, Corollary 2.2.14]:
$${\bold F_\bullet}:\; 0 \longrightarrow F_\ell \overset {\phi_\ell} \to {\longrightarrow }
F_{\ell-1} \longrightarrow \cdots \longrightarrow F_3 \overset {\phi_3}\to \longrightarrow F \overset \phi\to \longrightarrow G \overset \psi\to \longrightarrow H.$$ By Theorem 1.5, $(\bold F_\bullet)_\alpha$ is an exact complex of free $R_\alpha$-modules for almost all $\alpha.$ Hence $F_\alpha \overset {\phi_\alpha}\to \longrightarrow G_\alpha \overset{\psi_\alpha} \to \longrightarrow H_\alpha$ is exact for almost all $\alpha.$ \qed \enddemo \medskip
 
\heading 2. Specialization of arbitrary modules \endheading \smallskip

In this section we shall introduce the concept of a specialization of an  arbitrary finitely generated $R$-module. Let $L$ be an arbitrary finitely generated $R$-module.

\proclaim{Definition}{\rm  Let $F_1 \overset \phi\to \longrightarrow F_0 \longrightarrow L \longrightarrow 0$ be a finite  free presentation of $L$. Let $\phi_\alpha:(F_1)_\alpha \longrightarrow (F_0)_\alpha$ be a specialization of $\phi$ (as defined in Section 1). We call $L_\alpha :=  \Coker \phi_\alpha$ a {\it specialization} of $L$ (with respect to $\phi).$ \endproclaim

If $E$ is a free $R$-module, then $0 \longrightarrow E \longrightarrow E \longrightarrow 0$ is a finite free presentation of $E.$ Hence $E_\alpha$ as defined in Section 1 is a specialization of $E.$ \smallskip

\demo{Remark}{\rm The above definition of specialization covers Seidenberg's 
definition for the specialization of submodules of free modules. Let $L$ be a 
submodule of a free $R$-module $E = R^s$. For every element 
$l(u,x) = \big(f_1(u,x) ,\ldots,f_s (u,x)\big) \in k[u,x]^s$ set 
$l(\alpha,x) := \big(f_1(\alpha,x), \ldots ,f_s (\alpha,x)\big).$ 
Seidenberg [S, Appendix] defined the specialization of $L$ as the submodule  
$L_\alpha^*$ of $E_\alpha = R_\alpha^s$ generated by the elements $l(\alpha,x)$, 
$l(u,x) \in L \cap k[u,x]^s$. Clearly, this  definition extends the 
specialization of an ideal (see Introduction). 
Let $l_1(u,x),\ldots,l_t(u,x) \in k[u,x]$ be a generating set of $L$ and 
$F_1 \overset \phi \to \longrightarrow F_0 \longrightarrow L \longrightarrow 0$
a finite free presentation of $L$ corresponding to this generating set. 
Assume that  $$l_i(u,x) := \big(f_{i1} (u,x),\ldots,f_{is} (u,x)\big),\; i= 1,\ldots,t.$$ Then we have an exact sequence $F_1 \overset {\phi}\to \longrightarrow F_0 \overset {\psi}\to \longrightarrow E$ where $\psi$ is given by the matrix $\big(f_{ij} (u,x)\big)^T$. By Corollary 1.6, the sequence $(F_1)_\alpha \overset {\phi_\alpha}\to \longrightarrow (F_0)_\alpha \overset {\psi_\alpha} \to \longrightarrow E_\alpha$ is exact for almost all $ \alpha $. Hence $\Im \psi_\alpha $ is the submodule of $E_\alpha$ generated by the elements $l_1(\alpha,x), \ldots,l_t(\alpha,x)$. By [S, Appendix, Theorem 1], the latter module is equal to $L_\alpha^*$ for almost all $\alpha $. Thus,
$$L_\alpha^* \cong \Im \psi_\alpha \cong \Coker \phi_\alpha = L_\alpha$$ for almost all $\alpha.$ \enddemo 
\smallskip

The definition of $L_\alpha$ clearly depends on the chosen presentations of $L.$ If we choose a different finite free presentation $F'_1 \longrightarrow F'_0  \longrightarrow L  \longrightarrow 0$ we may get a  different specialization $L'_\alpha$ of $L.$  We shall see that $L_\alpha$ and $L'_\alpha$ are canonically isomorphic. For this we need to introduce the specialization of a homomorphism of finitely generated $R$-modules. \smallskip 
Let $v:  L \longrightarrow M$ be a homomorphism of finitely generated $R$-modules. Consider a commutative diagram $$\CD F_1 @ >\phi>>F_0@ >>>L @>>> 0 \\ @VVv_1V   @VVv_0V  @VVvV\\ G_1 @ >\psi>>G_0@ >>>M @>>> 0 \endCD $$ where the rows are finite free presentations of $L, M$. For almost all $\alpha$ we have the diagram $$\CD (F_1)_\alpha @>\phi_\alpha>> (F_0)_\alpha@>>>L_\alpha @>>> 0\\ @VV(v_1)_\alpha V   @VV(v_0)_\alpha V  @VVv_\alpha V \\
(G_1)_\alpha @ >\psi_\alpha>>(G_0)_\alpha@ >>>M_\alpha @>>> 0 \endCD$$ where $(v_0)_\alpha  \phi_\alpha = \psi_\alpha  (v_1)_\alpha$ by Lemma 1.2 and $v_\alpha$ is the induced homomorphism which makes the above diagram commutative. 

\proclaim{Definition} {\rm The induced homomorphism $v_\alpha\,:\, L_\alpha \longrightarrow M_\alpha$ is called a {\it specialization} of $v.$}\endproclaim
This definition depends only on the presentations of $L$ and $M$.

\proclaim{Lemma 2.1}  $v_\alpha$ does not depend on the choices of $v_0$ and $v_1$ for almost all $\alpha$.\endproclaim 
\demo{Proof} Suppose that we are given two other maps $w_i : F_i \longrightarrow G_i$, $i =0, 1, $ lifting the same homomorphism $v : L  \longrightarrow M.$ By [E, A3.13], there exists a homomorphism $h : F_0 \longrightarrow G_1$ such that $v_0 - w_0 = \psi h.$ By Lemma 1.2, 
$(v_0)_\alpha-(w_0)_\alpha = \psi_\alpha h_\alpha$ for almost all $\alpha.$ Hence the maps $(w_i)_\alpha$ induce the same map $v_\alpha \,:\, L_\alpha \longrightarrow M_\alpha.$ \qed\enddemo \smallskip
>From now on we will identify the given $R$-modules $L, M$ with fixed finite free presentations 
$F_1 \longrightarrow F_0 \longrightarrow L \longrightarrow 0,\; G_1 \longrightarrow G_0 \longrightarrow M \longrightarrow 0.$ A homomorphism $v: L \longrightarrow M$ will correspond to a map between these finite presentations (as complexes) $$\CD F_1 @ >>>F_0 @>>>L @>>> 0\\ @VVv_1V @VVv_0V  @VVvV\\ G_1 @ >>>G_0 @ >>>M @>>> 0. \endCD$$ Here we will use the same symbol $v$ with subindex $i$ to denote the given map from $F_i$ to $G_i,$ $i = 0, 1.$ We shall see that there are canonical isomorphisms between the specializations of an $R$-module with respect to different finite free presentations. 

\proclaim {Proposition 2.2} Let $v: L \longrightarrow M$ be an isomorphism of finitely generated $R$-modules. Then, for almost all $\alpha,$ the specialization $v_\alpha: L_\alpha \longrightarrow M_\alpha$  is an isomorphism of $R_\alpha$-modules.\endproclaim 
\demo {Proof} Put $w = v^{-1}$. There is a commutative diagram of the form $$\CD F_1 @ >>>F_0@ >>>L @>>> 0 \\ @VV v_1V @VV v_0V  @VVvV\\ G_1 @ >>>G_0@ >>>M @>>> 0\\ @VV w_1V  @V Vw_0V  @VVwV \\ F_1 @ >>>F_0 @ >>>L@ >>> 0.\\ \endCD $$ The composed map $\{w_i v_i | \; i= 0, 1 \}$ of complexes of free $R$-modules gives a presentation of the identity map $\text {id}_L: L \longrightarrow L.$ Hence it is homotopy equivalent to the map $\{\text {id}_{F_i}| \; i = 0, 1\}$, by [E, A 3.14]. By Lemma 1.2, the map $\{(w_i)_\alpha (v_i)_\alpha | \; i= 0, 1 \}$ of complexes of free $R_\alpha$-modules is homotopy equivalent to the map $\{\text {id}_{(F_i)_\alpha}| \; i = 0, 1\}.$ Thus, they lift to the same map from $L_\alpha$ to $L_\alpha,$ that is $v_\alpha w_\alpha = \text{id}_{L_\alpha}$ for almost all $\alpha.$ Similarly, we also have $w_\alpha v_\alpha = \text {id}_{L_\alpha}$ for almost all $\alpha.$  \qed \enddemo \smallskip
Let $L_\alpha$ and $L'_\alpha$ be specializations of $L$ with respect to different finite free presentations. By Proposition 2.2, $L_\alpha \cong L'_\alpha$ and the isomorphism is just the specialization of the identity map $\text{id}_L: L \longrightarrow L$ with respect to the two given presentations of $L.$ In particular, the specialization of the identity map $\text {id}_L: L \longrightarrow L$ with respect to a fixed finite free presentation of $L$ is the identity map $\text{id}_{L_\alpha}: L_\alpha \longrightarrow L_\alpha.$ \smallskip
If we choose two different finite free presentations for the map $v: L \longrightarrow M$ we may get two different specializations $v_\alpha: L_\alpha \longrightarrow M_\alpha$, $v'_\alpha: L'_\alpha \longrightarrow M'_\alpha$ of $v.$ However, there are isomorphisms $\varepsilon_L$ and $\varepsilon_M$ such that the diagram $$\CD L_\alpha @>v_\alpha >>M_\alpha\\ @ V\varepsilon_LVV @V\varepsilon_MVV\\ L'_\alpha @>v'_\alpha >>M'_\alpha\\ \endCD$$ is commutative. In fact, $\varepsilon_L$ and $\varepsilon_M$ are the specializations of the identity maps $\text{id}_L$ and $\text{id}_M$. The commutativity of the diagram follows from the second statement of the following result.

\proclaim {Lemma 2.3} Let $v, w: L \longrightarrow M,\,z : M \longrightarrow N$ be homomorphisms of finitely generated $R$-modules. Then, for almost all $\alpha,$ \smallskip
\noindent\rm (i)\quad $(v + w)_\alpha = v_\alpha + w_\alpha,$ \smallskip
\noindent\rm (ii)\quad $(zv)_\alpha = z_\alpha v_\alpha.$ \endproclaim
\demo{Proof} (i) From the maps between the finite free presentations of $L$ and $M$ corresponding to $v$ and $w$ we get the commutative diagram 
$$\CD F_1 @ >>>F_0 @ >>>L @>>> 0 \\ @VVv_1+w_1V @VVv_0+w_0V @VVv+wV \\ G_1 @ >>>G_0 @ >>>M @>>> 0. \\ \endCD $$
By Lemma 1.2, $(v_i + w_i)_\alpha = (v_i)_\alpha + (w_i)_\alpha$ for almost all $\alpha$, $i = 0,1$. Hence $(v+w)_\alpha = v_\alpha + w_\alpha$ for almost all $\alpha.$ \smallskip
\noindent (ii) Consider a commutative diagram
$$\CD F_1 @ >>>F_0@ >>>L @> >> 0 \\ @VV v_1V @VV v_0V @VVvV \\ G_1 @ >>>G_0@ >>>M @ >>> 0 \\ @VV z_1V  @VVz_0V  @VVzV \\ H_1 @ >>>H_0 @ >>>N @ >>> 0, \\ \endCD $$ 
where the last rows is a finite free presentation of $N.$ By Lemma 1.2, $(z_i v_i)_\alpha = (z_i)_\alpha (v_i)_\alpha$ for almost all $\alpha$, $ i = 0, 1.$ Hence $(zv)_\alpha = z_\alpha v_\alpha$ for almost all $\alpha.$ \qed\enddemo \smallskip

We now come to the main result of this section.
\proclaim {Theorem 2.4} Let $0\longrightarrow L \longrightarrow M \longrightarrow N \longrightarrow 0 $ be an exact sequence of finitely generated $R$-modules. Then the sequence $0\longrightarrow L_\alpha \longrightarrow M_\alpha \longrightarrow N_\alpha\longrightarrow 0$
is exact for almost all $\alpha.$ \endproclaim 
\demo{Proof} By [E, Corollary 19.8] or [BH, Corollary 2.2.14], given two arbitrary finite free presentations $F_1 \longrightarrow F_0 \longrightarrow L \longrightarrow 0$ and $H_1 \longrightarrow H_0 \longrightarrow N \longrightarrow 0$, we can extend them to finite free resolutions $\bold F_\bullet$ and $\bold H_\bullet$ of $L$ and $N$. Then $\bold F_\bullet \oplus \bold H_\bullet$ is a finite free resolution of $M.$ By Lemma 1.2, the following diagram is commutative for almost all $\alpha$:
$$\CD @. 0 @. 0 @. @. 0 \\ @.@VVV @VVV @. @VVV \\ 0@>>>(F_l)_\alpha @ >(\phi_l)_\alpha >> (F_{l-1})_\alpha@ >>>\ldots @>(\phi_1)_\alpha >>(F_0)_\alpha \\ @. @VVV @VVV @. @VVV  \\ 0 @>>>(F_l\oplus H_l)_\alpha @ >(\psi_l)_\alpha >>(F_{l-1}\oplus H_{l-1})_\alpha @ >>>\ldots @>(\psi_1)_\alpha >>(F_0\oplus H_0)_\alpha \\ @. @VVV    @VVV @. @VVV \\ 0@>>>(H_l)_\alpha @ >(\delta_l)_\alpha >>(H_{l-1})_\alpha @ >>>\ldots @>(\delta_1)_\alpha >>(H_0)_\alpha \\ @. @VVV @VVV @. @VVV \\ @. 0 @. 0 @. @. 0 \\ \endCD $$
where the rows and columns are exact for almost all $\alpha $ by Theorem 1.5. Using induction we can easily prove that the sequence
$$0 \longrightarrow \text{\rm Coker}\,(\phi_i)_\alpha \longrightarrow \text{\rm Coker}\,(\psi_i)_\alpha \longrightarrow \text{\rm Coker}\,(\delta_i)_\alpha \longrightarrow 0 $$ is exact for $i = l,l-1,\ldots ,1.$ For $i = 1$ we have $\Coker (\phi_1)_\alpha = L_\alpha,\; \Coker (\delta_1)_\alpha = N_\alpha,$ and $\text{\rm Coker}\,(\psi_1)_\alpha$ is the specialisation of $M$ with respect to the finite free presentation $F_1 \oplus H_1 \longrightarrow F_0 \oplus H_0 \longrightarrow M \longrightarrow 0.$ Thus, there is a commutative diagram 
$$\CD
0 @ >>>L_\alpha @ >>>\text{\rm Coker}(\psi_1)_\alpha @>>>N_\alpha @>>>0 \\
@. \| @.  @VV\varepsilon_\alpha V \| \\ 0 @ >>>L_\alpha @ >>>M_\alpha @>>>N_\alpha @>>>0 \\ \endCD $$
where $M_\alpha$ is the specialization of $M$ with respect to an arbitrary finite presentation and $\varepsilon_\alpha$ the specialization of the identity map $\text{id}_M: M \longrightarrow M$. By Proposition 2.2, $\varepsilon_\alpha$ is an isomorphism. By Lemma 2.3, the above diagram is commutative. Therefore, as the first sequence is exact, so is the second sequence. \qed  \enddemo \smallskip 

\proclaim {Corollary 2.5} Let $v : L \longrightarrow M$ be a homomorphism of finitely generated $R$-modules. Then, for almost all $\alpha,$ \smallskip  
\noindent {\rm (i)}\quad $\Ker v_\alpha \cong (\Ker v)_\alpha$ \smallskip
\noindent {\rm (ii)} \quad $\Im v_\alpha \cong (\Im v)_\alpha $ \smallskip
\noindent {\rm (iii)} \quad $\Coker v_\alpha \cong (\Coker v)_\alpha.$ \endproclaim
\demo {Proof} By Theorem 2.4, the sequences $$0\longrightarrow (\Ker v)_\alpha \longrightarrow L_\alpha \longrightarrow (\Im v)_\alpha \longrightarrow 0 $$
$$0\longrightarrow (\Im v)_\alpha \longrightarrow M_\alpha \longrightarrow (\Coker v)_\alpha \longrightarrow 0 $$ are exact for almost all $\alpha.$ Now consider the diagram: 
$$\matrix 0 \longrightarrow (\Ker v)_\alpha \longrightarrow  L_\alpha & & \overset v_\alpha \to \longrightarrow & & M_\alpha  \longrightarrow (\Coker v)_\alpha \longrightarrow 0 \\ & \searrow & & \nearrow & \\ & &  (\Im v)_\alpha & & \endmatrix$$ By Lemma 2.3 (ii), the triangle involving $v_\alpha$ is commutative for almost all $\alpha.$ Hence the conclusions are immediate. \qed \enddemo \smallskip

\demo{Remark} By Corollary 2.5, if $v$ is injective (resp. surjective) then $v_\alpha$ is injective (resp. surjective) for almost all $\alpha.$\enddemo 

\proclaim {Corollary 2.6} Let $L$ be a finitely generated $R$-module. Then 
$$\proj L_\alpha = \proj L$$ for almost all $\alpha.$  \endproclaim
 
\demo{Proof} If $\text{proj.dim}\,L = 0,$ then $L$ is a free $R$-module. Hence $L_\alpha$ is a free $R_\alpha$-module by definition and $\text{proj.dim}\,L_\alpha = 0.$ If $\text{proj.dim}\,L > 0,$ $L$ is not a free $R$--module. By [B-H, Proposition 1.4.9], $I(\phi) \ne R.$ Let $F_1 \overset \phi\to \longrightarrow F_0 \longrightarrow L \longrightarrow 0$ be a finite free presentation of $L.$ By [S, Appendix, Theorem 3], $I(\phi_\alpha) \ne R_\alpha$ for almost all $\alpha.$ Hence $L_\alpha$ is not a free $R_\alpha$--module. Put $H = \text {Im}\,\phi.$ Then $\text{proj.dim}\,H = \text{proj.dim}\,L - 1.$ Using induction we may assume that proj.dim\,$H_\alpha $ = proj.dim\, $H$ for almost all $\alpha.$  By Theorem 2.4, the sequence $0 \longrightarrow H_\alpha \longrightarrow (F_0)_\alpha \longrightarrow L_\alpha \longrightarrow 0$ is exact for almost all $\alpha $. Hence $$\proj L_\alpha  = \proj H_\alpha + 1  = \proj H +1 = \proj L$$ for almost all $\alpha.$ \qed \enddemo \smallskip

\proclaim {Proposition 2.7} Let ${\bold F_\bullet}:\; 0 \longrightarrow F_\ell \overset {\phi_\ell} \to \longrightarrow  F_{\ell-1} \longrightarrow \cdots \longrightarrow F_1 \overset {\phi_1} \to \longrightarrow F_0 \longrightarrow 0$ be a finite complex of finitely generated $R$-modules. Then, for almost all $\alpha,$ $$({\bold F_\bullet})_\alpha :\; 0 \longrightarrow  (F_\ell)_\alpha \overset (\phi_\ell)_\alpha \to \longrightarrow (F_{\ell-1})_\alpha \longrightarrow \cdots \longrightarrow (F_1)_\alpha \overset {(\phi_1)_\alpha}\to \longrightarrow (F_0)_\alpha \longrightarrow 0 $$ is a finite complex of $R_\alpha$-modules with $[H_i({\bold F_\bullet})]_\alpha = H_i(({\bold F_\bullet})_\alpha).$ \endproclaim

\demo{Proof} By Corollary 2.5, since the map $\Im \phi_i \longrightarrow \Ker \phi_{i-1}$ is injective, the map $(\Im \phi_i)_\alpha \longrightarrow (\Ker \phi_{i-1})_\alpha$ and therefore the map $\Im (\phi_i)_\alpha \longrightarrow \Ker (\phi_{i-1})_\alpha$ is injective, too. So $({\bold F_\bullet})_\alpha$ is a finite complex of $R_\alpha$-modules. By Theorem 2.4, 
$$[H_i({\bold F_\bullet})]_\alpha = (\Ker \phi_{i-1} / \Im \phi_i)_\alpha = (\Ker \phi_{i-1})_\alpha / (\Im \phi_i)_\alpha.$$ Using Corollary 2.5 again we obtain
$$(\Ker \phi_{i-1})_\alpha / (\Im \phi_i)_\alpha \cong \Ker (\phi_{i-1})_\alpha / \Im (\phi_i)_\alpha \cong H_i(({\bold F_\bullet})_\alpha). \qed $$ \enddemo

\heading 3. Operations on specializations of modules \endheading \smallskip  

So far we have defined the specialization of an $R$-module, but we have not discussed what operations are allowed between them. In this section we will discuss some of the basic operations on specializations of modules. 

\proclaim {Lemma 3.1} Let $L$ and $M$ be finitely generated $R$-modules. There is a finite free presentation of $L \oplus M$ such that for almost all $\alpha$, $(L \oplus M)_\alpha = L_\alpha \oplus M_\alpha.$\endproclaim
\demo{Proof} Let $F_1 \overset {\phi} \to \longrightarrow F_0 \longrightarrow L \longrightarrow 0$ and $G_1 \overset {\psi} \to \longrightarrow G_0 \longrightarrow M \longrightarrow 0$ be free presentations of $L$ and $M.$ Then $$F_1 \oplus G_1 \overset {\phi \oplus \psi} \to \longrightarrow F_0 \oplus G_0 \longrightarrow L \oplus M \longrightarrow 0$$ is a free presentation of $L \oplus M.$ By definition and Lemma 1.3, $$(L \oplus M)_\alpha = \Coker\,\phi_\alpha \oplus \psi_\alpha = \Coker\,\phi_\alpha \oplus \Coker\,\psi_\alpha = L_\alpha \oplus M_\alpha.\qed$$ \enddemo 

Let us now consider a submodule $M$ of a finitely generated $R$-module $L$. We shall see that $M_\alpha$  may be identified with a submodule of $L_\alpha.$ Indeed, the canonical map $M_\alpha \longrightarrow L_\alpha$ is injective for almost all $\alpha$ by Corollary 2.5. Moreover, if we fix a finite free presentation of $L$, then different finite free presentations of $M$ yield different specializations $M_\alpha, M'_\alpha$ with the same image in $L_\alpha.$ This follows from the commutative  diagram:
$$\matrix M_\alpha & & \\ & \searrow & \\ \wr \Vert \varepsilon_\alpha & & L_\alpha \\
& \nearrow & \\ M'_\alpha & & \endmatrix $$ where $\varepsilon_\alpha$ is the specialization of the identity map $\text{id}_M: M \longrightarrow M$ with respect to the two different presentations of $M.$ \smallskip  

By the above observation we may consider operations on specializations of submodules of a module.
\proclaim{Proposition 3.2} Let $M, N$ be submodules of a finitely generated $R$-module $L$. For almost all $\alpha$, there are canonical isomorphisms: \smallskip
\noindent \rm (i)\quad $(L/M)_\alpha \cong L_\alpha/M_\alpha$, \smallskip
\noindent \rm (ii)\quad $(M \cap N)_\alpha \cong M_\alpha \cap N_\alpha$, \smallskip
\noindent \rm (iii)\quad $(M+N)_\alpha \cong M_\alpha + N_\alpha$.\endproclaim 
\demo{Proof} (i) By Theorem 2.4 we have, for almost all $\alpha,$ the exact sequence 
$$0 \longrightarrow M_\alpha \longrightarrow L_\alpha \longrightarrow (L/M)_\alpha \longrightarrow 0,$$ 
\noindent which implies that there is a canonical isomorphism $(L/M)_\alpha \cong L_\alpha/M_\alpha.$ \smallskip
\noindent (ii) Consider the exact sequence $$ 0 \longrightarrow M \cap N \longrightarrow L \longrightarrow (L/M)\oplus (L/N) \longrightarrow 0.$$  \noindent By Lemma 3.1, $[(L/M)\oplus (L/N)]_\alpha = (L_\alpha/M_\alpha) \oplus (L_\alpha /N_\alpha).$ Hence we have, for almost all $\alpha,$ the exact sequence
$$ 0 \longrightarrow (M \cap N)_\alpha \longrightarrow L_\alpha \longrightarrow (L_\alpha /M_\alpha) \oplus (L_\alpha/N_\alpha) \longrightarrow 0$$ which implies that there is a canonical isomorphism $(M \cap N)_\alpha \cong M_\alpha \cap N_\alpha.$ \smallskip
\noindent (iii) Consider the exact sequence $$ 0 \longrightarrow M \oplus N \overset {\varphi} \to \longrightarrow L \longrightarrow L/M + N \longrightarrow 0$$ where $\varphi(x,y) = x+y, \;x \in M,\; y \in N.$ Similarly as above we can deduce that there are canonical isomorphisms  $(L/M + N)_\alpha \cong L_\alpha/M_\alpha + N_\alpha$ and, therefore, $(M+N)_\alpha \cong M_\alpha + N_\alpha$ for almost all $\alpha.$ \qed \enddemo 
As a consequence we shall see that generators of an $R$-module specializes.
\proclaim{Corollary 3.3} Let $L = Re_1 + \cdots + Re_s$ be an $R$-module. Then, for almost all $\alpha,$ there exist elements $(e_1)_\alpha,\cdots,(e_s)_\alpha \in L_\alpha$ such that  $(Re_j)_\alpha = R_\alpha (e_j)_\alpha, \; j = 1,\cdots,s,$ and $L_\alpha = R_\alpha (e_1)_\alpha + \cdots +  R_\alpha (e_s)_\alpha.$ \endproclaim
\demo{Proof} For each $j = 1,\cdots,s$ the canonical map $R \longrightarrow Re_j$ is surjective. By Corollary 2.5, $R_\alpha \longrightarrow (Re_j)_\alpha$ is surjective for almost all $\alpha.$ Hence we can find $(e_j)_\alpha \in (Re_j)_\alpha \subseteq L_\alpha$ such that $(Re_j)_\alpha = R_\alpha(e_j)_\alpha.$ Now, applying Proposition 3.2 (iii) we obtain $$L_\alpha = (Re_1)_\alpha + \cdots +  (Re_s)_\alpha = R_\alpha (e_1)_\alpha + \cdots +  R_\alpha (e_s)_\alpha. \qed$$ \enddemo

\demo{Remark} If $L = I = (f_1(u,x),\cdots,f_s(u,x))$ is an ideal of $R$ and $e_j = f_j(u,x),$ then $(e_j)_\alpha = f_j(\alpha,x),\; j = 1,\cdots,s.$ Thus, $I_\alpha = (f_1(\alpha,x),\cdots,f_s(\alpha,x))$ as already shown by Krull [Kr1, \S1, Satz 1]. \enddemo \smallskip

An important property of specializations of modules is that they preserve the dimension.

\proclaim{Theorem 3.4} Let $L$ be a finitely generated $R$-module. Then, for almost all $\alpha,$ \smallskip
\noindent \rm (i)\quad $\Ann L_\alpha = (\Ann L)_\alpha,$ \smallskip
\noindent \rm (ii)\quad $\dim L_\alpha = \dim L.$
\endproclaim
\demo{Proof} (i) Let $L = Re_1 + \cdots + Re_s.$ Then $Re_j = R/\text{\rm Ann}\; e_j ,\; j = 1,\cdots,s.$ By Theorem 2.4, we have the exact sequence $$0 \longrightarrow (\Ann e_j)_\alpha \longrightarrow R_\alpha \longrightarrow (Re_j)_\alpha \longrightarrow 0$$ for almost all $\alpha.$ By Corollary 3.3, there exists $(e_j)_\alpha \in L_\alpha$ such that $(Re_j)_\alpha = R_\alpha(e_j)_\alpha.$ It follows that $\Ann (e_j)_\alpha = (\Ann e_j)_\alpha.$ Since  $L_\alpha = R_\alpha (e_1)_\alpha + \cdots +  R_\alpha (e_s)_\alpha,$ we have $$\Ann L_\alpha = \bigcap\limits_{j = 1}^s \Ann (e_j)_\alpha = \bigcap\limits_{j = 1}^s (\Ann e_j)_\alpha =(\bigcap\limits_{j = 1}^s \Ann e_j)_\alpha = (\Ann L)_\alpha.$$
\noindent (ii) We have $\height \Ann L_\alpha = \height (\Ann L)_\alpha = \height \Ann L$ by [T, Lemma 1.1]. Therefore $\dim L_\alpha = \dim L$ for almost all $\alpha.$ \qed \enddemo \smallskip

\proclaim{Proposition 3.5} Let $k$ be an algebraically closed field. Let $M$ be a submodule of a finitely generated $R$-module $L.$ Let $M = \bigcap\limits_{j = 1}^s M_j$ be a primary decomposition, where $M_j$ are primary submodules of $L$ with associated prime ideals $\wp_j, \; j = 1, \ldots,s.$ For almost all $\alpha,$ $M_\alpha = \bigcap\limits_{j = 1}^s (M_j)_\alpha $ is a primary decomposition, where $(M_j)_\alpha$ are primary submodules of $L_\alpha$ with associated prime ideals $(\wp_j)_\alpha, \; j = 1, \ldots,s.$ \endproclaim
\demo{Proof} By Proposition 3.2 (ii), $M_\alpha = \bigcap\limits_{j = 1}^s (M_j)_\alpha$ for almost all $\alpha.$ Since $k$ is algebraically closed, the residue field of $R/\wp_j$ is a regular extension of $k(u).$ By [T, Proof of Lemma 2.1], this implies that $(\wp_j)_\alpha$ is a prime ideal for almost all $\alpha.$ Since $\wp_j = \sqrt{\Ann (L/M_j)},$
$$(\wp_j)_\alpha = \big(\sqrt{\Ann (L/M_j)}\ \big)_\alpha = \sqrt{\Ann (L/M_j)_\alpha} = \sqrt{\Ann (L_\alpha/(M_j)_\alpha)}$$ by Theorem 3.4 (i) and Proposition 3.2 (i). Therefore, we can conclude that $(M_j)_\alpha$ is a $(\wp_j)_\alpha$-primary submodule of $L_\alpha$ for almost all $\alpha.$ \qed \enddemo \smallskip
\demo {Remark} If $k$ is not algebraically closed, Proposition 3.5 does not hold. For example, if $I$ is the prime ideal $(x^2 + 1)$ of ${\Bbb R}(u)[x],$ then $I_\alpha = (x+i)(x-i)$ for all $\alpha$ such that ${\Bbb R}(\alpha) = \Bbb C.$ \enddemo \smallskip 

\proclaim{Proposition 3.6} Let $L$ be a finitely generated $R$-module and $I$ an ideal of $R$. For almost all $\alpha$ we have
$$\align  (0_L : I)_\alpha &\cong 0_{L_\alpha} : I_\alpha,\\ (IL)_\alpha  &\cong I_\alpha L_\alpha .\endalign $$ \endproclaim

\demo{Proof} We first prove the statements for the case $I$ is a principal ideal $(f(u,x))$. Consider the commutative diagram 
$$\CD F_1 @>>>F_0 @>>>L @>>>0\\ @VVf(u,x)V  @VVf(u,x)V    @VVf(u,x) V\\ F_1 @>>>F_0 @>>>L @>>>0.\\ \endCD $$
\noindent By definition, the specialization of the map $F_i \overset {f(u,x)}\to \longrightarrow F_i$ is the map $(F_i)_\alpha \overset {f(\alpha,x)}\to \longrightarrow (F_i)_\alpha.$ Hence the specialization of the map $L \overset {f(u,x)}\to \longrightarrow L$  is the map $L_\alpha \overset {f(\alpha,x)}\to \longrightarrow L_\alpha.$ By Theorem 2.4, we obtain the exact sequence 
$$\CD 0 @>>> (0_L : f(u,x))_\alpha @>>> L_\alpha @> {f(\alpha,x)}>> L_\alpha \endCD$$
for almost all $ \alpha.$ Thus, $$(0_L : f(u,x))_\alpha \cong 0_{L_\alpha}: f(\alpha,x),$$ and therefore
$$[f(u,x)\,L]_\alpha \cong [L/0_L : f(u,x)]_\alpha \cong L_\alpha /0_{L_\alpha} : f(\alpha,x) \cong f(\alpha,x)\, L_\alpha$$ for almost all $\alpha.$ \par
Now assume that $I = (f_1(u,x),\ldots,f_s(u,x))$. Then $I_\alpha = (f_1(\alpha,x),\ldots,f_s(\alpha,x))$ for almost all $\alpha$ by Corollary 3.3. Since 
$$0_L : I = \bigcap\limits_{i = 1}^s 0_L : f_i(u,x),$$ using Proposition 3.2 (ii) we obtain
$$\align (0_L : I)_\alpha &= \big(\bigcap\limits_{i = 1}^s 0_L : f_i (u,x)\big)_\alpha = \bigcap\limits_{i = 1}^s (0_L : f_i (u,x))_\alpha \cr &\cong \bigcap\limits_{i = 1}^s 0_{L_\alpha} : f_i (\alpha,x)  = 0_{L_\alpha} : I_\alpha \endalign$$ 
for almost all $\alpha.$ Similarly, we have 
$$ \align (IL)_\alpha &= (\sum\limits^s_{i=1} f_i (u,x)L)_\alpha = \sum\limits^s_{i=1} [f_i (u,x)L]_\alpha \\ &\cong \sum\limits^s_{i=1} f_i (\alpha ,x)L_\alpha = I_\alpha L_\alpha \endalign$$ for almost all $\alpha.$ \qed \enddemo \medskip

\heading 4. Specialization of Tor and Ext \endheading \smallskip

First we will show that for any free $R$-module $E$ of finite rank and any finitely generated $R$-module $M,$ there is a canonical map $(E\otimes_R M)_\alpha \longrightarrow E_\alpha \otimes_{R_\alpha} M_\alpha$ for almost all $\alpha.$ Note that similarly, there is a canonical map $\Hom_R(E,M)_\alpha \longrightarrow \Hom_{R_\alpha}(E_\alpha,M_\alpha)$ for almost all $\alpha.$

\proclaim{ Lemma 4.1} Let $\phi : E \longrightarrow F$ be a homomorphism of free $R$-modules of finite ranks and $M$ a finitely generated $R$-module. Then, for almost all $\alpha,$ there are isomorphisms 
$$\align \varepsilon_E : & (E\otimes_R M)_\alpha \longrightarrow E_\alpha\otimes_{R_\alpha} M_\alpha \\ \varepsilon_F : & (F\otimes_R M)_\alpha \longrightarrow F_\alpha\otimes_{R_\alpha} M_\alpha \endalign $$ 
such that the following diagram is commutative:
$$\CD (E\otimes_R M)_\alpha @>(\phi \otimes \text{id}_M)_\alpha >>(F \otimes_R M)_\alpha \\ @V\varepsilon_EVV  @V\varepsilon_FVV\\ E_\alpha \otimes_{R_\alpha} M_\alpha @>\phi_\alpha \otimes \text{id}_{M_\alpha} >>F_\alpha\otimes_{R_\alpha} M_\alpha  \\ \endCD $$ \endproclaim 
\demo {Proof} Let $G_1 \longrightarrow G_0 \longrightarrow M \longrightarrow 0$ be a finite free presentation $M.$ Then the sequence $E {\otimes}_R G_1 \longrightarrow E {\otimes}_R G_0 \longrightarrow E {\otimes}_R M \longrightarrow 0$ is a finite free presentation of $E {\otimes}_R M.$ By Lemma 1.3, we get the commutative diagram
$$\CD (E \otimes_R G_1)_\alpha @>>> (E \otimes_R G_0)_\alpha @>>> (E \otimes_R M)_\alpha @>>> 0 \\ @| @| @. \\ E_\alpha \otimes_{R_\alpha} (G_1)_\alpha @>>> E_\alpha \otimes_{R_\alpha} (G_0)_\alpha @>>> E_\alpha \otimes_{R_\alpha} M_\alpha @>>> 0\\ 
\endCD $$
where the last row is induced by the exact sequence $(G_1)_\alpha \longrightarrow (G_0)_\alpha \longrightarrow M_\alpha \longrightarrow 0.$ The above diagram induces an isomorphism $ \varepsilon_E: (E {\otimes}_R M)_\alpha \longrightarrow E_\alpha {\otimes}_{R_\alpha} M_\alpha $ for almost all $\alpha.$ Similarly, there is an isomorphism $\varepsilon_F: (F {\otimes}_R M)_\alpha \longrightarrow F_\alpha {\otimes}_{R_\alpha} M_\alpha $ for almost all $\alpha.$ From the commutative diagrams \smallskip
$$\CD E_\alpha \otimes_{R_\alpha}(G_1)_\alpha  @>>> E_\alpha \otimes_{R_\alpha}(G_0)_\alpha @>>> E_\alpha \otimes_{R_\alpha}M_\alpha @>>> 0\\ @| @| @V(\varepsilon_E)^{-1}VV @. \\ (E \otimes_R G_1)_\alpha  @>>> (E \otimes_R G_0)_\alpha @>>> (E \otimes_R M)_\alpha @>>> 0\\ @VVV  @VVV @V(\phi \otimes \text{id}_M)_\alpha VV @. \\ (F \otimes_R G_1)_\alpha  @>>> (F \otimes_R G_0)_\alpha @>>> (F \otimes_R M)_\alpha @>>> 0\\ @| @| @V\varepsilon_FVV @. \\ F_\alpha \otimes_{R_\alpha}(G_1)_\alpha  @>>> F_\alpha \otimes_{R_\alpha}(G_0)_\alpha @>>> F_\alpha \otimes_{R_\alpha}M_\alpha @>>> 0 \endCD $$ and  
$$\CD E_\alpha \otimes_{R_\alpha}(G_1)_\alpha  @>>> E_\alpha \otimes_{R_\alpha}(G_0)_\alpha @>>> E_\alpha \otimes_{R_\alpha}M_\alpha @>>> 0\\
@VVV  @VVV @V\phi_\alpha \otimes \text{id}_{M_\alpha}VV @. \\
 F_\alpha \otimes_{R_\alpha}(G_1)_\alpha  @>>> F_\alpha \otimes_{R_\alpha}(G_0)_\alpha @>>> F_\alpha \otimes_{R_\alpha}M_\alpha @>>> 0 \endCD $$
we get 
$$\varepsilon_F (\phi \otimes \text{id}_M)_\alpha (\varepsilon_E)^{-1} = \phi_\alpha \otimes \text{id}_{M_\alpha}$$ for almost all $\alpha$. \qed\enddemo \smallskip

\proclaim{Theorem 4.2} Let $L, M$ be finitely generated $R$-modules. Then, for almost all $\alpha,$ $$\Tor_i^R(L,M)_\alpha \cong \Tor_i^{R_\alpha}(L_\alpha,M_\alpha),\; i \ge 0.$$ \endproclaim
\demo {Proof}Let ${\bold F_\bullet}: \; 0 \longrightarrow F_\ell \longrightarrow 
F_{\ell-1} \longrightarrow \cdots \longrightarrow F_1 \longrightarrow F_0$ be a free resolution of $L.$ By Proposition 2.7, 
$$\Tor_i^R(L,M)_\alpha = [H_i({\bold F_\bullet} \otimes_R M)]_\alpha 
\cong H_i(({\bold F_\bullet} \otimes_R M)_\alpha)$$ for almost all $\alpha.$ By Theorem 1.5, 
$$({\bold F_\bullet})_\alpha : \; 0 \longrightarrow (F_\ell)_\alpha  \longrightarrow  (F_{\ell-1})_\alpha \longrightarrow \cdots \longrightarrow (F_1)_\alpha  \longrightarrow (F_0)_\alpha $$ is a free resolution of $L_\alpha$ for almost all $\alpha.$ Therefore 
$$\Tor_i^{R_\alpha}(L_\alpha,M_\alpha) = H_i(({\bold F_\bullet})_\alpha \otimes_{R_\alpha} M_\alpha)$$
for almost all $\alpha.$ By Lemma 4.1, 
$$({\bold F_\bullet} \otimes_R M)_\alpha \cong ({\bold F_\bullet})_\alpha \otimes_{R_\alpha} M_\alpha$$
for almost all $\alpha.$ Hence we can conclude that $\Tor_i^R(L,M)_\alpha \cong \Tor_i^{R_\alpha}(L_\alpha,M_\alpha)$ for almost all $\alpha.$\qed \enddemo \smallskip

\proclaim{ Theorem 4.3} Let $L,M$ be finitely generated $R$-modules. Then, for almost all $\alpha,$
$$\Ext_R^i(L,M)_\alpha \cong \Ext_{R_\alpha}^i(L_\alpha,M_\alpha),\; i\ge 0.$$ \endproclaim
\demo{Proof} Take a free resolution ${\bold F_\bullet}$ of $L.$ By Proposition 2.7, 
$$ \Ext_R^i(L,M)_\alpha = [H^i(\text{Hom}_R({\bold F_\bullet},M))]_\alpha \cong H^i(\Hom_R({\bold F_\bullet},M)_\alpha)$$ 
for almost all $\alpha.$ By Theorem 1.5, $({\bold F_\bullet})_\alpha$ is a free resolution of $L_\alpha$ for almost all $\alpha.$ Therefore, 
$$\Ext_{R_\alpha}^i(L_\alpha,M_\alpha) = H^i(\Hom_R(({\bold F_\bullet})_\alpha,M_\alpha).$$
Similarly as in the proof of Lemma 4.1 we can prove that 
$$\Hom_R({\bold F_\bullet},M)_\alpha \cong \Hom_{R_\alpha}(({\bold F_\bullet})_\alpha,M_\alpha)$$ 
for almost all $\alpha.$ Hence $\Ext_R^i(L,M)_\alpha \cong  \Ext_{R_\alpha}^i(L_\alpha ,M_\alpha)$ for almost all $\alpha.\qed $ \enddemo \smallskip

Let $L$ be a finitely generated $R$-module and $I$ an ideal of $R$. The {\it grade} of $I$ on $L$, denoted by $\gra (I,L)$, is the maximal length of regular $L$-sequences in $I$ if $L \ne IL$. If $L = IL,$ $\gra (I,L) := \infty$. In particular, the {\it grade} of $L$ is the grade of $\Ann L$ on $R,$ denoted by $\gra L.$ We will show that the grade remains unchanged by specialization.

\proclaim{ Corollary 4.4}  Let $L$ be a finitely generated $R$-module. Then, for almost all $\alpha,$ \smallskip
\noindent {\rm (i)}\quad $\gra (I_\alpha,L_\alpha) = \gra (I,L)$ \smallskip
\noindent {\rm (ii)}\quad $\gra L_\alpha = \gra L$ for any ideal $I$ of $R.$
\endproclaim
\demo{Proof }(i) Without restriction we may assume that $L \ne IL.$ By [B-H, Theorem 1.2.5], 
$$\align \gra (I,L)&= \min \{i\; |\;\Ext_R^i(R/I,L) \ne 0 \},\\
\gra (I_\alpha,L_\alpha) &= \min \{i\; |\;\Ext_{R_\alpha}^i(R_\alpha /I_\alpha,L_\alpha) \ne 0 \}.\endalign$$ By Theorem 4.3, we have 
$$\min \{i\; |\;\Ext_R^i(R/I,L) \ne 0 \} = \min \{i\; |\;\Ext_{R_\alpha}^i(R_\alpha /I_\alpha,L_\alpha) \ne 0 \}$$ for almost all $\alpha.$ Hence the conclusion is immediate. \par
\noindent (ii) Using Theorem 3.4 we obtain 
$$\align \gra L_\alpha &= \gra (\Ann L_\alpha,R_\alpha) = \gra ((\Ann L)_\alpha,R_\alpha) \\ &= \gra (\Ann L,R) = \gra L.\qed \endalign$$ \enddemo 
Recall that a finitely generated $R$-module $L$ is {\it perfect} if $\gra L = \proj L.$ It is well-known that perfect $R$-modules are Cohen-Macaulay modules.\smallskip
\proclaim{ Proposition 4.5} Let $L$ be a finitely generated perfect $R$-module. Then $L_\alpha$ is a perfect $R_\alpha$-module for almost all $\alpha.$\endproclaim
\demo {Proof} This is a consequence of Corollary 4.4 (ii) and Corollary 2.6. \qed \enddemo \smallskip
\proclaim{Corollary 4.6} Let $L$ be a finitely generated Cohen-Macaulay $R$-module. Assume that $\text{\rm supp}\, L$ is connected. Then $L_\alpha$ is a Cohen-Macaulay $R_\alpha$-module for almost all $\alpha.$\endproclaim
\demo {Proof} By [B-V, Proposition (16.19)], the assumption implies that $L$ is perfect. Thus, $L_\alpha$ is perfect and hence Cohen-Macaulay for almost all $\alpha.$ \qed \enddemo  \smallskip

\heading Acknowledgement \endheading
 The authors are grateful to the referee for his suggestions. \vglue 0.8cm

\heading References \endheading \smallskip

{\parindent=0.95cm 
\item{\hbox to 0.85cm{[B-H]\hfill}} W. Bruns and J. Herzog, {\it Cohen-Macaulay rings,} Cambridge University Press, 1993. \smallskip
\item{\hbox to 0.85cm{[B-V]\hfill}} W. Bruns and U. Vetter, {\it Determinantal rings}, Lect. Notes Math. 1327, Springer, 1988.\smallskip 
\item{\hbox to 0.85cm{[B-E]\hfill}} D.A. Buchsbaum and D. Eisenbud, {\it What makes a complex exact?} J. Algebra 25 (1973), 259-268. \smallskip
\item{\hbox to 0.85cm{[E]\hfill}}  D. Eisenbud, {\it Commutative algebra with a view toward algebraic geometry, } Springer, 1995. \smallskip
\item{\hbox to 0.85cm{[Kr1]\hfill}} W. Krull, {\it Parameterspezialisierung in Polynomringen,} Arch. Math.~1 (1948), 56-64. \smallskip
\item{\hbox to 0.85cm{[Kr2]\hfill}} W. Krull, {\it Parameterspezialisierung in Polynomringen {\rm II}, Grundpolynom}, Arch. Math. 1 (1948), 129-137. \smallskip
\item{\hbox to 0.85cm{[Ku1]\hfill}} W.E. Kuan, {\it A note on a generic hyperplane section of an algebraic variety,} Can. J. Math. 22 (1970), 1047-1054. \smallskip
\item{\hbox to 0.85cm{[Ku2]\hfill}} W.E. Kuan, {\it Specialization of a generic hyperplane section through a rational point of an algebraic variety,} Ann. Math. Pura Appl. 94 (1972), 75-82.   \smallskip
\item{\hbox to 0.85cm{[S]\hfill}}  A. Seidenberg, {\it The hyperplane sections of normal varieties, } Trans. Amer. Math. Soc. 69 (1950), 375-386. \smallskip
\item{\hbox to 0.85cm{[T]\hfill}} N. V. Trung, {\it Spezialisierungen  allgemeiner Hyperfl\"achenschnitte und Anwendungen}, in: Seminar D.Eisenbud/B.Singh/W.Vogel, Vol.~1, Teubner-Texte zur Mathematik, Band 29 (1980), 4-43.
\enddocument